\documentclass[12pt]{article}
\usepackage{hyperref}
\usepackage{amsmath}
\usepackage{ntheorem}
\usepackage{mathtools}
\usepackage{epsf}
\usepackage{theorem}
\usepackage{indentfirst}
\usepackage{latexsym}
\usepackage{amssymb}
\usepackage[dvips]{graphicx}
\usepackage{flafter}
\usepackage{caption}
\usepackage{textcomp}
\usepackage{supertabular}
\usepackage{longtable, booktabs}
\usepackage{hhline}
\usepackage{bigstrut,bigdelim}
\usepackage{rotating}
\usepackage{multicol}
\usepackage{multirow}
\usepackage{makeidx}
\usepackage{epsfig}
\usepackage{fancyhdr}
\newif\ifger
\gerfalse
\newtheorem{theorem}{Theorem}
\newtheorem{lemma}{Lemma}[section]
\newtheorem{corollary}{Corollary}[section]

\newtheorem*{remark}{Remark}

\newtheorem{proposition}{Proposition}[section]
\newtheorem{example}{Example}[section]

\begin{document}
\baselineskip=19pt

\title{Block-transitive automorphism groups on 3-designs with small block size}

\author{Xiaoqin Zhan$^\dag$, Meifang Yang\footnote{Corresponding author. E-mail:  mfyang0918@126.com(M. Yang), zhanxiaoqinshuai@126.com(X. Zhan)}\\
\small \it  $\dag$ School of Science, East China JiaoTong University, \\
\small \it  Nanchang, 330013, PR China\\
\small \it $\ast $ School of Statistics, Jiangxi University of Finance and Econmics,\\
\small \it Nanchang, 330013, PR China}

\date{}
\maketitle

\begin{abstract}
The paper is an investigation of the structure of block-transitive automorphism groups of a 3-design with  small block size. Let $G$  be a block-transitive automorphism group of a nontrivial $3$-$(v,k,\lambda)$ design $\mathcal{D}$ with  $k\le 6$. We prove that if $G$ is point-primitive then $G$ is of affine or almost simple type. If $G$ is point-imprimitive then $\mathcal{D}$ is a  $3$-$(16,6,\lambda)$ design with $\lambda\in\{4, 12, 16, 24, 28, 48, 56, 64, 84, 96, 112, 140\}$, and $rank(G)=3$.

\medskip
\noindent{\bf MSC(2010):} Primary: 05B25; Secondary: 20B25.

\medskip
\noindent{\bf Keywords:} Block-transitive; Automorphism group; 3-design; Point-primitive; Point-imprimitive

\end{abstract}

\section{Introduction}
A $t$-$(v, k, \lambda)$ {\it design} is a pair $(\mathcal{P}, \mathcal{B})$ in which $\mathcal{P}$ is a $v$-set of points and $\mathcal{B}$ is a collection of $k$-sets of $\mathcal{P}$ called {\it blocks}, such that every $t$-set of $\mathcal{P}$
is contained in precisely $\lambda$ blocks. If $t<k<v-1$ holds, then we speak of a {\it nontrivial} $t$-design. It is {\it simple} if no two blocks are identical. All of the $t$-designs in this paper will be simple and nontrivial.

An {\it automorphism} of ${\mathcal D}$ is a permutation of $\cal P$ which leaves $\mathcal B$ invariant. The {\it full automorphism group}  of $\mathcal{D}$ consists of all automorphisms of $\mathcal{D}$ and is denoted by ${\rm Aut}(\mathcal{D})$. A subgroup $G$ of the automorphism group of $\mathcal D$ is  {\it block-transitive} if it acts transitively  on
$\mathcal{B}$; $\mathcal D$ is said to be block-transitive if Aut($\mathcal D$) is. {\it Point-} and {\it flag-transitivity} are defined similarly. A set of blocks of  $\mathcal{D}$ is called a set of {\it base blocks} with respect to an automorphism group $G$ of $\mathcal{D}$ if it contains exactly one block from each $G$-orbit on the block set. In particular, if $G$ is a block-transitive automorphism group of $\mathcal{D}$, then any block $B$ is a base block of $\mathcal{D}$.

Block-transitivity is just one of many conditions that can be imposed on the automorphism group $G$
of a $t$-design $\cal D$. It is well known that if $G$ is block-transitive, then $G$ is also point-transitive (Block's Lemma \cite{Block}). It is elementary that the flag-transitivity of $G$ on  a linear space ($2$-$(v,k,1)$ design) implies its point-primitivity.   By a result of Davies \cite{Dav1}, for  $2$-$(v,k,\lambda)$ designs, this implication remains true if $(r,\lambda)=1$ (where $r$ denotes the number of blocks containing a given point). However, block-transitivity does not necessarily imply  point-primitivity. For example, let $\cal D$ be a 2-design consisting of the points and hyperplanes of any Desarguesian projective space $PG(n,q)$ where $n\geq2$ and $(q^{n+1}-1)/(q-1)$ is not a prime, and take $G$ as the group generated by a Singer cycle.

If the autommorphism group $G$ of $\cal D$ is point-primitive, then $G$ is of one of the following five types by O'Nan-Scott theorem (see \cite{Lie1988} for details).
\begin{enumerate}
  \item[\rm (i)]  Affine.
  \item[\rm (ii)]  Almost simple.
    \item[\rm (iii)]  Product.
  \item[\rm (iv)] Simple diagonal.
  \item[\rm (v)] Twisted wreath product.

\end{enumerate}

 In 1984, Camina and Gagen \cite{CG} proved that if $G$ is block-transitive on a $2$-$(v,k,1)$ design $\cal{D}$ with $k\mid v$, then $G$ is either  point-primitive of affine or almost simple type. Inspired by the proof, several others \cite{BDD,Dav2,Zies1988} generalised the result in \cite{CG} to prove that groups acting flag-transitively on $2$-$(v,k,1)$ designs are affine or almost simple. It is worth nothing that both \cite{Dav2} and \cite{Zies1988} generalised the result to the situation of 2-designs  with $(r,\lambda)=1$.
For a $t$-$(v,k,\lambda)$ design, Cameron and Praeger \cite{Cameron} proved the following result in 1993.
 \begin{proposition}  \label{CP}
Let $\mathcal{D}=(\mathcal{P}, \mathcal{B})$ be a   $t$-$(v,k,\lambda)$  design with $t\ge2$. Then the following holds:
\begin{enumerate}
\item[\rm(i)] \,  If $G\le {\rm Aut}(\cal D)$ acts block-transitively on $\mathcal{D}$, then $G$ also acts $\lfloor t/2\rfloor$-homogeneously on $\cal P$.
\item[\rm(ii)] \, If $G\le {\rm Aut}(\cal D)$ acts flag-transitively on $\mathcal{D}$, then $G$ also acts  $\lfloor (t+1)/2\rfloor$-homogeneously on $\cal P$.
 \end{enumerate}
\end{proposition}
According to this result, if $G$ acts block-transitively on a $t$-$(v,k,\lambda)$ design $\cal D$ with $t\geq4$ then $G$ is either point-primitive of affine or almost simple type as $G$ is 2-homogeneous on the points of $\cal D$. Therefore, it is necessary to study the block-transitive $t$-$(v,k,\lambda)$ designs with $t\le3$.

 The main aim of this paper is to study $3$-$(v,k,\lambda)$ designs admitting a block-transitive automorphism group $G$.
Firstly, we analyse the case in which the automorphism group $G$ is point-primitive, and we prove a reduction theorem for small values of $k$.
\begin{theorem}\label{th1}
Let $G$ be a block-transitive automorphism group of a nontrivial $3$-$(v, k, \lambda)$ design with $k\leq6$. If  $G$ is point-primitive, then $G$ is of affine type, or almost simple type.

\end{theorem}

In fact, there exist many $3$-designs admitting a block-transitive, point-primitive automorphism group of affine or almost simple type. Here are some examples (cf. \cite{Huber}):
\begin{example} {\rm
\begin{itemize}
  \item[\rm (i)] Let $\mathcal{D}=(\mathcal{P},\mathcal{B})$, where  $\cal P$ and  $\cal B$ are the points  and planes of the affine space $AG(d,2)$ with $d\ge3$. Then $\cal D$ is isomorphic to the $3$-$(2^d,4,1)$ design admitting $G=AGL(d,2)$ as its flag-transitive (block-transitive), point-primitive automorphism group of affine type.
  \item[\rm (ii)] Let $\mathcal{D}$ be the Mathieu-Witt $3$-$(22,6,1)$ design, and $G\unrhd M_{22}$. Then $G$ is a flag-transitive (block-transitive), point-primitive automorphism group of $\cal D$ with almost simple action.
\end{itemize}}
\end{example}

For the point-imprimitive case, Delandtsheer and Doyen have shown in \cite{DD} that if  $\cal D$ is a  $t$-$(v,k,\lambda)$ design admitting a block-transitive point-imprimitive automorphism group $G$ then $v\le(\binom{k}{2}-1)^2$. Assume that $G$ has a system of $d$ blocks of imprimitivity each of size $c$.
 In \cite[Corollaries 3.2 and 3.4]{Cameron}, it was shown that for a  block-transitive, point-imprimitive $3$-$(v,k,\lambda)$ design with $d=2$ or $c=2$ then $v\leq\binom{k}{2}+1$. Thus, for a fixed block size $k$, there are only finitely many $t$-$(v,k,\lambda)$ designs with a block-transitive automorphism group which is point-imprimitive.

Secondly, the other purpose of this paper is to study  $3$-$(v,k,\lambda)$ designs  admitting a  block-transitive point-imprimitive automorphism group and prove the following theorem:

\begin{theorem}\label{th2}
 Let $\mathcal{D}=(\cal P, \cal B)$ be a  nontrivial $3$-$(v,k,\lambda)$ design with $k\le 6$ and admitting a block-transitive automorphism group $G$. If $G$ is point-imprimitive then $rank(G)=3$, and $\mathcal{D}$ is a $3$-$(16,6,\lambda)$ design  with $$\lambda\in\{4,12,16,24,28,48,56,64,84,96,112,140\}.$$
\end{theorem}

The paper is organized as follows. In Section 2, we introduce some preliminary results that are important for the remainder of the paper. In Sections 3 and 4, we shall give the proofs of the Theorems \ref{th1} and \ref{th2} respectively.


\section{Preliminaries}

 The notation and terminology used is standard and can be found in  \cite{Handbook,Demb1968} for design theory and in \cite{Dixon,Wie1964} for group theory. In particular, if $G$ is a permutation group on point set $\cal P$, and $\alpha\in B\subseteq \mathcal{P}$, then $G_{\alpha}$ denotes the stabilizer of a point $\alpha$ in $G$, and $G_B$ denotes the setwise stabilizer of $B$ in $G$, and
 $G_{\alpha B}$ denotes the stabilizer of a flag $(\alpha,B)$ in $G$.

\begin{lemma}{\rm \cite[1.2, 1.9]{Handbook}} \label{L1}
The parameters $v, b,  r, k, \lambda$ of a $3$-design
satisfy the following conditions:
\begin{enumerate}
\item[\rm(i)] \, $vr=bk$.
\item[\rm(ii)] \, $\lambda v(v-1)(v-2)=bk(k-1)(k-2)$.
\end{enumerate}
\end{lemma}

\bigskip

  The following lemma is useful for the study of block-transitive $3$-$(v,k,\lambda)$ designs.

\begin{lemma} \label{divide}
 Let $\mathcal{D}=(\cal P, \cal B)$ be a  nontrivial $3$-$(v,k,\lambda)$ design with $k\le 6$ and admitting a block-transitive automorphism group $G$. Then $r$ divides  $k|G_{\alpha}|$. Furthermore, $r$ divides  $k\lambda d(d-1)$, and  $(v-1)(v-2)$ divides $k(k-1)(k-2)d(d-1)$, for all nontrivial subdegrees $d$ of $G$.
 \end{lemma}
\textbf{Proof.}\quad  Let $B$ be a block of $\mathcal{D}$ containing the point $\alpha$. The point-transitivity and block-transitivity imply
$$|G:G_{\alpha B}|=|G:G_\alpha||G_\alpha:G_{\alpha B}|=v|G_\alpha:G_{\alpha B}|,$$
and
$$|G:G_{\alpha B}|=|G:G_B||G_B:G_{\alpha B}|=b|G_B:G_{\alpha B}|.$$
 Hence, $|G_\alpha:G_{\alpha B}|=\frac{r|G_B:G_{\alpha B}|}{k}$ by Lemma \ref{L1}(i), and so $r$ divides $k|G_{\alpha}|$.
In order to prove the remaining result, here we  prove only the case  $k=4$, and  the result in lemma can be proved imitate to the proof of the case $k=4$ for the other values of $k$.

Clearly, $|G_B:G_{\alpha B}|=1,2,3$ or 4 as $k=|B|=4$. We will analyze each of these cases separately.

$\mathbf{(1)}$ Let $|G_B:G_{\alpha B}|=1$, then $|G_\alpha:G_{\alpha B}|=\frac{r}{4}$.  Suppose that $G_\alpha$ has four orbits with same size on pencil $P(\alpha)$ (i.e. blocks containing a given point $\alpha$), and denoted by $\mathcal{O}_1$, $\mathcal{O}_2$, $\mathcal{O}_3$ and $\mathcal{O}_4$, respectively.
Let $\Gamma\neq\{\alpha\}$ be a nontrivial $G_\alpha$-orbit with $|\Gamma|=d$. Set $\mu_i=|\Gamma\cap B_i|$ where $B_i\in \mathcal{O}_i$ ($i=1,2,3,4$). Clearly, $0\leq\mu_i\leq3$. Counting the number of set $\{(\{\beta,\gamma\}, B)\mid \{\beta,\gamma\}\in B\cap \Gamma\}$ in two ways, and we get
$$\frac{r}{4}\sum_{i=1}^4\binom{\mu_i}{2}=\lambda\binom{d}{2}.$$
So $r$ divides $4\lambda d(d-1).$

Suppose that $G_\alpha$ has three orbits $\mathcal{O}_1$, $\mathcal{O}_2$, and $\mathcal{O}_3$ with sizes $\frac{r}{4},\frac{r}{4}$ and $\frac{r}{2}$ on pencil $P(\alpha)$ respectively. Then
\begin{equation}\label{O3}
\frac{r}{4}\binom{\mu_1}{2}+\frac{r}{4}\binom{\mu_2}{2}+\frac{r}{2}\binom{\mu_3}{2}=
\lambda\binom{d}{2}.
\end{equation}
Also, $r$ divides $4\lambda d(d-1).$

 Assume that $G_\alpha$ has two orbits $\mathcal{O}_1$ and $\mathcal{O}_2$ with $|\mathcal{O}_1|=\frac{r}{4}$ and $|\mathcal{O}_2|=\frac{3r}{4}$ on pencil $P(\alpha)$, we obtain
  \begin{equation}\label{O2}
  \frac{r}{4}\binom{\mu_1}{2}+\frac{3r}{4}\binom{\mu_2}{2}=\lambda\binom{d}{2}.
  \end{equation}
  Hence, $r$ divides $4\lambda d(d-1).$

$\mathbf{(2)}$ Let  $|G_B:G_{\alpha B}|=2$, then $|G_\alpha:G_{\alpha B}|=\frac{r}{2}$. If $G_\alpha$ has three orbits with sizes $\frac{r}{2}$, $\frac{r}{4}$, $\frac{r}{4}$ then the Equation (\ref{O3}) holds. If $G_\alpha$ has two orbits $\mathcal{O}_1$ and $\mathcal{O}_2$ with sizes $\frac{r}{2}$ and $\frac{r}{2}$ respectively, then $$\frac{r}{2}\binom{\mu_1}{2}+\frac{r}{2}\binom{\mu_2}{2}=\lambda\binom{d}{2}.$$
In both of cases we have $r$ divides $4\lambda d(d-1)$.

$\mathbf{(3)}$ Let  $|G_B:G_{\alpha B}|=3$, then $|G_\alpha:G_{\alpha B}|=\frac{3r}{4}$, and $G_\alpha$ has two orbits with sizes $\frac{3r}{4}$ and $\frac{r}{4}$ respectively. Thus, $r$ divides $4\lambda d(d-1)$ by Equation (\ref{O2}).

$\mathbf{(4)}$ Let  $|G_B:G_{\alpha B}|=4$, then $|G_\alpha:G_{\alpha B}|=r$ and so $G_{\alpha}$ acts transitively on $P(\alpha)$. Set $\mu=|\Gamma\cap B|$ where $B\in P(\alpha)$. We obtain $$r\binom{\mu}{2}=\lambda\binom{d}{2}.$$
Hence $r$ divides $\lambda d(d-1)$, and so  $4\lambda d(d-1)$ is divisible by $r$.

By Lemma \ref{L1}(i)(ii), $r=\frac{\lambda (v-1)(v-2)}{(k-1)(k-2)}$ and so $(v-1)(v-2)$ divides $24d(d-1)$.
 $\hfill\square$

 \bigskip

From the proof of  \cite[Lemma 2.3]{Alavi} we get the following:
\begin{lemma}\label{neq}
There does not exist  a non-abelian finite simple group satisfying
$$(|T|-1)(|T|-2)<480|{\rm Out}(T)|.$$
\end{lemma}

\bigskip

In the study of point-imprimitive case, the basis of our method is the following elementary result.
\begin{lemma}{\rm \cite[Proposition 1.1]{Cameron}}\label{ele}
Let $\mathcal{D}=(\mathcal{P},\mathcal{B})$ be a $t$-$(v,k,\lambda)$ design, admitting a block-transitive automorphism group $G$. Let $H$ be a permutation group with $G\le H\le S_v$, and $\mathcal{B}^*=\mathcal{B}^H$ the set of images of blocks in $\cal B$ under $H$. Then $(\mathcal{P},\mathcal{B}^*)$ is a $t$-$(v,k,\lambda^*)$ design, for some $\lambda^*$, admitting the block-transitive automorphism group $H$.
\end{lemma}

\section{Primitivity}
 The principal tool used in the proof is the O'Nan-Scott theorem for finite primitive groups proved by Liebeck, Praeger and Saxl in \cite{Lie1988}.  We will prove Theorem \ref{th1} by dealing with
the cases of  product action, simple diagonal action and twisted wreath product action
 separately. The proof of the Theorem \ref{th1} is inspired by the proof  of \cite[Theorem 1.1]{TianZhou2013}.

\subsection{Product action}
Here, we suppose that $G$ has a product action on $\mathcal{P}$. Then  $G\leq K^{m}\rtimes S_{m}=K\wr  S_{m}$ with  $m\geq2$, where $K$ is a primitive group (of almost simple or
diagonal type) on $\Omega$ of size $v_{0}\geq5$, and $\mathcal{P}=\Omega^{m}$.

\begin{proposition}\label{pa}  Let $\mathcal{D}=(\cal P, \cal B)$ be a  nontrivial $3$-$(v,k,\lambda)$ design with $k\leq6$ admitting a block-transitive point-primitive automorphism group $G$. Then $G$ is not of product action type.
\end{proposition}

{\bf Proof.}\, Assume the contrary, suppose that $H=K\wr  S_{m}$ with $S_{m}$ acting on the set $M=\{1,2,\ldots,m\}$.  Let $\alpha$ and $\beta$ be two distinct points of $\cal P$. Then $d=|\beta^{G_{\alpha}}|$ is a subdegree of $G$. Since $G$ is a subgroup of $H$, it follows that
\begin{equation}\label{gongshi2}
d=|G_{\alpha}:G_{\alpha\beta}|\leq|H_{\alpha}:H_{\alpha\beta}|.
\end{equation}
Let $\alpha=(\gamma,\gamma,\ldots,\gamma)\in\mathcal{P}$,
$\beta=(\delta,\gamma,\ldots,\gamma)\in\mathcal{P}$ with $\delta\neq\gamma$
and let $B\cong K^{m}$ be the base group of $H$. Then
$B_{\alpha}=K_{\gamma}^m$,
$B_{\alpha\beta}=K_{\gamma\delta}\times K_{\gamma}^{m-1}$. Now $H_{\alpha}=K_{\gamma}\wr S_{m}$,
and $H_{\alpha\beta}\geq K_{\gamma\delta}\times(K_{\gamma}\wr
S_{m-1})$. Suppose $K$ has rank $s$ on $\Omega$ with $s\geq2$. We can choose a  $\delta\in\Omega$ satisfying $|K_{\gamma}:K_{\gamma\delta}|\leq\frac{v_{0}-1}{s-1}$, so that
\begin{equation*}
|H_{\alpha}:H_{\alpha\beta}|=\frac{|H_{\alpha}|}
{|H_{\alpha\beta}|}
\leq\frac{|K_{\gamma}|^{m}\cdot m!}{|K_{\gamma\delta}||K_{\gamma}|^{m-1}\cdot(m-1)!}\leq m\frac{v_{0}-1}{s-1},
\end{equation*}
and hence  $d\leq m\frac{v_{0}-1}{s-1}$ by  Equation (\ref{gongshi2}). From Lemma \ref{divide} we have
\begin{equation*}\label{gongshi3}
  (v-1)(v-2)\le k(k-1)(k-2)\cdot m\frac{v_{0}-1}{s-1}\cdot(m\frac{v_{0}-1}{s-1}-1).
\end{equation*}
Combining this with $v=v_0^m$ and $k\leq6$ we get all possible $(v_0, m, s)$ as in Table \ref{TAB1}.
  \begin{table}[htbp]
 \centering
\caption{ All possible values of  $v_0, m, s$ with $k\leq6$}\label{TAB1}
\begin{tabular}{clll}
\toprule
&$k=4$ & $k=5$ & $k=6$ \\
\midrule
$(m,s)=(2,2)$&$ v_0 \in\{5,6,7,8\}$&$v_0\in\{5,6,\ldots,14\}$&$v_0\in\{5,6,\ldots,20\}$\\
$(m,s)=(2,3)$&$\emptyset$&$v_0\in\{5,6\}$&$ v_0\in\{5,6,7,8,9\}$\\
$(m,s)=(3,2)$&$\emptyset$&$\emptyset$&$ v_0=5$\\
\bottomrule
\end{tabular}
\end{table}

First, assume that $(m,s)=(2,2)$.  Then $K$ acts 2-transitively on $\Omega$, and $H=K\wr S_{2}$ has rank 3 with subdegrees 1, $2(v_{0}-1)$, $(v_{0}-1)^{2}$ on the point set $\mathcal{P}=\Omega\times\Omega$. Note that $G\leq H$, so each subdegree of $H$ is the sum of some subdegrees of $G$,  by Lemma \ref{divide} we conclude that $(v_0^2-1)(v_0^2-2)$ divides $k(k-1)(k-2)(2v_0-2)(2v_0-3)$, it is impossible.

For the case  $(m,s)=(2,3)$, $K$ is a primitive group with rank 3  on $\Omega$. From \cite[9.62 Table]{Handbook}, there is no such group $K$ with a primitive action (of almost simple or
diagonal type) and  rank 3 on a set $\Omega$ of size $ v_{0}\in\{5,6,7,8,9\}$.

Now,  assume that $(m,s)=(3,2)$. Then $H=K\wr S_{3}$ has rank 4 with subdegrees 1, $3(v_{0}-1)$, $3(v_{0}-1)^2$, $(v_{0}-1)^{3}$ on the point set $\mathcal{P}=\Omega\times\Omega\times\Omega$. This contradicts the fact that $(v_0^3-1)(v_0^3-2)$ divides $k(k-1)(k-2)(3v_0-3)(3v_0-4)$ as $v_0=5$ and $k=6$.  $\hfill\square$

\subsection{Simple diagonal action}

Suppose that $G$ is a primitive group of simple diagonal type. Then
$M ={\rm Soc}(G) = T_1\times\cdots\times T_m \cong T^m$ and $M_{\alpha}\cong T$ is a diagonal subgroup of $M$, where $T_i \cong T$ is a non-abelian finite simple group, for $i = 1, \ldots, m$ and $m\geq2$.
Here  $G_{\alpha}$ is isomorphic to a subgroup of ${\rm Aut}(T) \times S_m$ and has an orbit $\Gamma$ in $\mathcal{P}-\{\alpha\}$ with $|\Gamma|\le m|T|$.

\begin{proposition} \label{sd}
  Let $\mathcal{D}=(\cal P, \cal B)$ be a nontrivial $3$-$(v,k,\lambda)$ design with $k\leq6$ admitting a block-transitive point-primitive automorphism group $G$. Then $G$ is not of simple diagonal type.
\end{proposition}
{\bf Proof.}\,  If $G$ is of simple diagonal type, then  $|\mathcal{P}| = |T|^{m-1}$ and $G$ has a subdegree $d$ less than $m|T|$. From Lemma \ref{divide}, we have $$(|T|^{m-1}-1)(|T|^{m-1}-2)\leq k(k-1)(k-2)\cdot m|T|\cdot(m|T|-1).$$
It is easy to get $m=2$ as $k\leq6$ and $|T|\geq60$.

Also by Lemma \ref{divide}, we have that $r$ divides $ k|G_\alpha|$, and so $k |{\rm Aut}(T)||S_2|$ is divisible by $r$. Since
 ${\rm Out}(T)\cong{\rm Aut}(T)/{\rm Inn}(T)$ and ${\rm Inn}(T)\cong T/Z(T)$, it yields that  $r$ divides $2k|T||{\rm Out}(T)|$ as $T$ is a non-abelian simple group. Combining this with Lemma \ref{L1}(ii), we get that
  $(|T|-1)(|T|-2)$ divides $2k(k-1)(k-2)|T||{\rm Out}(T)|$.
Then $(|T|, |T|-1)=1$ and $(|T|, |T|-2)=2$ imply $$(|T|-1)(|T|-2)<4k(k-1)(k-2)
|{\rm Out}(T)|\leq480|{\rm Out}(T)|.$$
This violates Lemma \ref{neq}. $\hfill\square$

\subsection{Twisted wreath product action}
Next, we suppose that $G$ is a primitive group of twisted wreath
product type on $\cal P$. Let $\alpha\in \mathcal{P}$. Then $G\cong \prescript{}{Q}{B}\rtimes P$, where $P=G_{\alpha}$ is a transitive permutation group on $\{1,\ldots,m\}$ with $m\geq6$,  and
$\prescript{}{Q}{B}={\rm Soc}(G)=T_{1}\times\cdots\times T_{m}\cong
T^{m}$ is regular for some nonabelian simple groups $T$. Thus, $v=|\mathcal{P}|=|T|^{m}$.
Moreover, $G_{\alpha}$ has an orbit $\Gamma$ with $|\Gamma|\leq m|T|$.
\begin{proposition} \label{tw}
 Let $\mathcal{D}=(\cal P, \cal B)$ be a nontrivial  $3$-$(v,k,\lambda)$ design with $k\leq6$ admitting a block-transitive point-primitive automorphism group $G$. Then $G$ is not of twisted wreath
product type.
\end{proposition}
{\bf Proof.}\, If $G$ is of twisted wreath product type, then the argument here is similar to the proof of Proposition \ref{sd}. By Lemma \ref{divide}, we easily observe that $$(|T|^{m-1}-1)(|T|^{m-1}-2)\leq k(k-1)(k-2)\cdot m|T|(m|T|-1).$$
Then the inequalities $k\leq6$ and $|T|\geq60$ imply $m\leq2$, this contradicts the fact that $m\geq 6$. $\hfill\square$

\bigskip

$\textbf{Proof of Theorem \ref{th1}}$ \ It follows from Propositions \ref{pa}-\ref{tw}.

\section{Imprimitivity}

 Suppose that $G$ is an  imprimitive group on the point set $\cal P$. Then $\cal P$ can be partitioned into  $d$ nontrivial blocks of imprimitivity $\Delta_j$, $j=1,\ldots,d$, each of size $c$, and so $v=|\mathcal {P}|=cd$, with $c, d>1$.  Let $B$ be a $k$-set of $\cal P$, and let $\mathcal{B}^*=B^G$. Then the sizes of the intersections of each element of $\mathcal{B}^*$ with the imprimitivity classes determine a partition of $k$, say $\mathbf{x} = (x_1, x_2, \ldots , x_d )$ with $x_1\ge x_2\ge \ldots \ge x_d $ and $\sum\limits_{i=1}^d x_i = k$. Set $b_t=\sum\limits_{i=1}^d x_i(x_i-1)\cdots(x_i-t+1)$. Note that $b_1=k$. By \cite[Proposition 2.2]{Cameron}, the following lemma holds.

\begin{lemma} \label{b2}
Let $\mathcal{D}^*=(\mathcal{P},\mathcal{B}^*)$. Then
 \begin{itemize}
   \item [{\rm(i)}] $\mathcal{D}^*$ is a $2$-design if and only if
   $$b_2=\sum_{i=1}^d x_i(x_i-1)=\frac{k(k-1)(c-1)}{(v-1)}.$$
   \item [{\rm (ii)}] $\mathcal{D}^*$ is a $3$-design if and only if it is a $2$-design and $$b_3=\sum_{i=1}^d x_i(x_i-1)(x_i-2)=\frac{k(k-1)(k-2)(c-1)(c-2)}{(v-1)(v-2)}.$$
 \end{itemize}

\end{lemma}

\begin{proposition}\label{imp}
Let $\mathcal{D}=(\cal P, \cal B)$ be a nontrivial $3$-$(v,k,\lambda)$ design with $k\leq6$ admitting a block-transitive point-imprimitive automorphism group $G$. Then $k=6$ and $v=16$.
\end{proposition}
\textbf{Proof.}\quad Since $S_c\wr S_d$ is a point-imprimitive maximal subgroup of symmetric group $S_v$, we only need to consider the case that the group $G=S_c\wr S_d$ acts point-imprimitively on $\cal D$ by Lemma \ref{ele}.

 Suppose that the block size $k=4$. Then the partition $\mathbf{x}$ of $k$ is $(3,1,0,\ldots,0)$ as $\cal D $ is a 3-design.
%
%
If $\mathbf{x}=(3, 1,0,\ldots , 0 )$ then there is no such pair $(c,d)$ satisfying Lemma \ref{b2}(i) as $c>1$ and $d>1$.

Now assume that $k=5$. By Lemma   \ref{b2}(i),  the partitions $\mathbf{x}$ of $k$ and parameters $c,d$ are listed in Table \ref{TAB2}.
  \begin{table}[htbp]
 \centering
\caption{ The partitions of $k=5$ and parameters $c$, $d$.}\label{TAB2}
\begin{tabular}{ccccc}
\toprule
$\mathbf{x}$&$(4, 1,\ldots , 0 )$&$(3, 1, 1 )$&$(3,2)$&$(2,2,1,\ldots,0)$\\
$(c,d)$&$\emptyset$& (7,3)&(3,2)&(2,3),(4,4)\\

\bottomrule
\end{tabular}
\end{table}
 By using Lemma \ref{b2}(ii), we get that $\mathbf{x}=(3,2)$ and $(c,d)=(3,2)$. Then $v=cd=6$, contradicts the nontriviality of $\cal D$.

Finally, we assume that $k=6$. Similarly, the the case $(5, 1,\ldots , 0 )$ does not happen by Lemma   \ref{b2}(i).  Other partitions $\mathbf{x}$ of $k$ and parameters $c,d$ are listed in Table \ref{TAB3}.
  \begin{table}[htbp]
 \centering
\caption{ The partitions of $k=6$ and parameters $c$, $d$.}\label{TAB3}
\begin{tabular}{cccccc}
\toprule
$\mathbf{x}$&$(4, 1, 1,\ldots , 0 )$&$(4,2)$&$(3,3)$&$(3,2,1,\ldots,0)$&$(3,1,1,1\ldots,0)$\\
$(c,d)$&(3,2)&(8,2)&(3,2)&$\emptyset$&$(2,3),(4,4)$\\

\bottomrule
\end{tabular}
\end{table}
From Lemma \ref{b2}(ii) and the nontriviality of $\cal D$, we get that $\mathbf{x}=(4,2)$ and $(c,d)=(8,2)$ and so $v=16$.  $\hfill\square$

\bigskip

\begin{corollary}\label{co}
Let $\mathcal{D}=(\cal P, \cal B)$ be a  $3$-$(16,6,\lambda)$ design admitting $G$ as its block-transitive, point-imprimitive automorphism group. Then $rank(G)=3$ with subdegrees $1,7,8$, and  $$\lambda\in\{4, 12, 16, 24, 28, 48, 56, 64, 84, 96, 112, 140\}.$$
\end{corollary}
\textbf{Proof.}\quad By the proof of Proposition \ref{imp}, we have that $G\leq S_8\wr S_2$. From Lemma \ref{ele}, there exists a $3$-$(16,6,\lambda^*)$  design $\mathcal{D}^*$ admitting $S_8\wr S_2$ as a  block-transitive, point-imprimitive automorphism group. Let $\Delta_1$ and $\Delta_2$ be the blocks of imprimitivity of $\cal P$, and let $\alpha\in \Delta_1$. Clearly, $\Delta_1^{G_{\alpha}}=\Delta_1$, so $|\Delta_1|$ is sum  of some subdegrees $d$ of $G$.  On the other hand, it follows from Lemma \ref{divide} that $7$ divides $d(d-1)$, and then we easily observe that $G_{\alpha}$ has subdegrees $1,7,8$.

Let $B$ be any base block of $\mathcal{D}^*$. Since the partition of block size $k$ is $\mathbf{x}=(4,2)$, without loss of generality, we set $|B\cap \Delta_1|=4$ and $|B\cap \Delta_2|=2$. Let $\mathcal{D}^*=(\mathcal{P},\mathcal{B}^*)$. Then each block of $\mathcal{D}^*$ is the 6-set in $\mathcal{P}$ with  partition $\mathbf{x}=(4,2)$ as $S_8\wr S_2$ acts transitively on $\{\Delta_1,\Delta_2\}$ and 4-transitively on $\Delta_i$ $(i=1,2)$.  Thus, the number of blocks in $\mathcal{D}^*$ is
$$|\mathcal{B}^*|=\binom{2}{1}\binom{8}{4}\binom{8}{2}=3920.$$
We further obtain $\lambda^*=140$ by Lemma \ref{L1}(ii), and so $\lambda\leq140$ as $\mathcal{B}\subseteq\mathcal{B}^*$.

 By using the software package {\sc Magma}\cite{Magma}-command  {\tt TransitiveGroups(16)}, we know that there are 1954 transitive groups on $\mathcal{P} = \{1,2,3,\ldots,16\}$, exactly 22 of which are primitive. Here we only consider that $G$ is one of  the remaining 1932 imprimitive groups. Note that, if $B$ is a base block of $\mathcal{B}$, then $B\in \mathcal{B}^*$. A simple calculation by using command {\tt Design<3,16|$B^{G}$>}, we get that  $\lambda\in\{4, 12, 16, 24, 28, 48, 56, 64, 84, 96, 112, 140\}.$
  $\hfill\square$

\bigskip

$\textbf{Proof of Theorem \ref{th2}}$ \ It follows from Proposition \ref{imp} and Corollary \ref{co}.

 \bigskip
 \begin{remark}\quad
 {\rm Up to isomorphism, there are $28$ different block-transitive point-imprimitive nontrivial $3$-$(v,k,\lambda)$ designs with $k\le6$ by using command {\tt IsIsomorphic(D1,D2)} (see Table \ref{TAB4}).
  \begin{table}[htbp]
 \centering
\caption{ The number of pairwise non-isomorphic 3-designs}\label{TAB4}
\begin{tabular}{ccccccccccccc}
\toprule
$\lambda$&4&12&16&24&28&48&56&64&84&96&112&140\\
$n$&5&4&5&1&1&6&1&1&1&1&1&1\\

\bottomrule
\end{tabular}
\end{table}
The notation $``n"$ in Table \ref{TAB4} means that there are $n$ pairwise non-isomorphic block-transitive point-imprimitive $3$-$(16,6,\lambda)$ designs.
}
 \end{remark}

\bigskip
\section*{Acknowledgements}
The authors would like to thank anonymous referees for providing us helpful and constructive comments and suggestions. This work is supported by the National Natural Science Foundation of China (Grant Nos. 11801174 and 11961026).


\begin{thebibliography}{10}
\bibitem{Alavi} S.H. Alavi, A. Daneshkhah, N. Okhovat, On flag-transitive automorphism groups
of symmetric designs, Ars Math. Contemp. 17 (2019), 619-626.

\bibitem{Magma} W. Bosma, J. Cannon, C. Playoust, The Magma Algebra System I: The User Language. J. Symb. Comput., 1997.
\bibitem{Block} R.E. Block, On the orbits of collineation groups, Math. Z. 96 (1967), 33-49.

\bibitem{BDD} F. Buekenhout, A. Delandtsheer and J. Doyen, Finite linear spaces with flag-transitive group, J. Combin. Theory Ser. A 49 (1988), 268-293.

 \bibitem{Cameron} P.J. Cameron and C.E. Praeger, Block-transitive
designs I: point-imprimitive designs, Discrete Math. 118 (1993), 33-43.


\bibitem{CG} A.R. Camina and T.M. Gagen, Block-transitive automorphism groups, J. Algebra 86 (1984), 549-554.
\bibitem{Handbook}  C.J. Colbourn, J.H. Dinitz, The CRC Handbook of Combinatorial Designs. CRC Press, Boca Raton, FL, 2007.

\bibitem{Dav1} H. Davies,  Flag-transitivity and primitivity,
Discrete Math. 63 (1987), 91-93.
\bibitem{Dav2} H. Davies,  Automorphisms of designs, PhD Thesis,
University of East Anglia, 1987.


\bibitem{DD} A. Delandtsheer and J. Doyen, Most block-transitive $t$-designs are point-primitive, Geom. Dedicata, 29 (1989), 307-310.


\bibitem{Demb1968} P. Dembowski, Finite Geometries, Springer-Verlag, New York, 1968.

\bibitem{Dixon} J.D. Dixon, B. Mortimer, Permutation Groups,  Springer-Verlag, New York, 1996.

\bibitem{Huber} M. Huber, The classification of flag-transitive Steiner $3$-designs, Adv. Geom. 5 (2005), 195-221.


\bibitem{Lie1988} M.W. Liebeck, C.E. Praeger and J. Saxl,
 On the O'Nan-Scott theorem for finite primitive permutation groups, J. Aust. Math. Ser. A 44 (1988), 389-396.



\bibitem{Wie1964} H. Wielandt, Finite Permutation Groups, Academic Press, New York, 1964.

\bibitem{TianZhou2013} D.L. Tian, S.L. Zhou, Flag-transitive point-primitive symmetric $(v, k, \lambda)$ designs with $\lambda$ at most 100, J. Combin. Des. 21 (2013), 127-141.

\bibitem{Zies1988} P.H. Zieschang, Flag-transitive automorphism groups
of 2-designs with $(r,\lambda)=1$, J. Algebra 118 (1988), 265-275.


\end{thebibliography}
\end{document}